\documentclass[11pt,fullpage]{article}

\newcommand{\ol}{\setlength{\itemsep}{0pt.}\begin{enumerate}}
\newcommand{\eol}{\end{enumerate}\setlength{\itemsep}{-\parsep}}
\newcommand{\ignore}[1]{}
\title{\bf On multivariate Newton-like inequalities}

\author{Leonid Gurvits \thanks{%
{\tt gurvits@lanl.gov}. Los Alamos National Laboratory, 
Los Alamos, NM. } 
}
\begin{document}


\begin{titlepage}

\maketitle

\begin{abstract}
We study multivariate entire functions and polynomials with non-negative coefficients.
A class of {\bf Strongly Log-Concave} entire functions, generalizing {\it Minkowski} volume polynomials, is introduced:
an entire function $f$ in $m$ variables is called {\bf Strongly Log-Concave}
 if the function
$(\partial x_1)^{c_1}...(\partial x_m)^{c_m} f$ is either zero or\\
$\log((\partial x_1)^{c_1}...(\partial x_m)^{c_m} f)$ is concave on $R_{+}^{m}$. We start with yet another point of view (of {\it propagation}) on the standard univarite
(or homogeneous bivariate) {\bf Newton Inequalities}. We prove analogues of the {\bf Newton Inequalities} in the multivariate {\bf Strongly Log-Concave} case. One of the corollaries of our new  Newton-like inequalities is the fact that the support $supp(f)$ of a {\bf Strongly Log-Concave} entire function $f$ is
discretely convex ($D$-convex in our notation). The proofs are based on a natural convex relaxation of the derivatives
$Der_{f}(r_1,...,r_m)$ of $f$ at zero and on the lower bounds
on $Der_{f}(r_1,...,r_m)$, which generalize the
{\bf Van Der Waerden-Falikman-Egorychev} inequality for the permanent of doubly-stochastic matrices. A few open questions are posed in the final section.

\end{abstract} 
\end{titlepage}
\newpage
 
\newtheorem{THEOREM}{Theorem}[section]
\newenvironment{theorem}{\begin{THEOREM} \hspace{-.85em} {\bf :} 
}%
                        {\end{THEOREM}}
\newtheorem{LEMMA}[THEOREM]{Lemma}
\newenvironment{lemma}{\begin{LEMMA} \hspace{-.85em} {\bf :} }%
                      {\end{LEMMA}}
\newtheorem{COROLLARY}[THEOREM]{Corollary}
\newenvironment{corollary}{\begin{COROLLARY} \hspace{-.85em} {\bf 
:} }%
                          {\end{COROLLARY}}
\newtheorem{PROPOSITION}[THEOREM]{Proposition}
\newenvironment{proposition}{\begin{PROPOSITION} \hspace{-.85em} 
{\bf :} }%
                            {\end{PROPOSITION}}
\newtheorem{DEFINITION}[THEOREM]{Definition}
\newenvironment{definition}{\begin{DEFINITION} \hspace{-.85em} {\bf 
:} \rm}%
                            {\end{DEFINITION}}
\newtheorem{EXAMPLE}[THEOREM]{Example}
\newenvironment{example}{\begin{EXAMPLE} \hspace{-.85em} {\bf :} 
\rm}%
                            {\end{EXAMPLE}}
\newtheorem{CONJECTURE}[THEOREM]{Conjecture}
\newenvironment{conjecture}{\begin{CONJECTURE} \hspace{-.85em} 
{\bf :} \rm}%
                            {\end{CONJECTURE}}
\newtheorem{PROBLEM}[THEOREM]{Problem}
\newenvironment{problem}{\begin{PROBLEM} \hspace{-.85em} {\bf :} 
\rm}%
                            {\end{PROBLEM}}
\newtheorem{QUESTION}[THEOREM]{Question}
\newenvironment{question}{\begin{QUESTION} \hspace{-.85em} {\bf :} 
\rm}%
                            {\end{QUESTION}}
\newtheorem{REMARK}[THEOREM]{Remark}
\newenvironment{remark}{\begin{REMARK} \hspace{-.85em} {\bf :} 
\rm}%
                            {\end{REMARK}}
\newtheorem{FACT}[THEOREM]{Fact}
\newenvironment{fact}{\begin{FACT} \hspace{-.85em} {\bf :} 
\rm}%
		            {\end{FACT}}

 
\newcommand{\thm}{\begin{theorem}}
\newcommand{\lem}{\begin{lemma}}
\newcommand{\pro}{\begin{proposition}}
\newcommand{\dfn}{\begin{definition}}
\newcommand{\rem}{\begin{remark}}
\newcommand{\xam}{\begin{example}}
\newcommand{\cnj}{\begin{conjecture}}
\newcommand{\prb}{\begin{problem}}
\newcommand{\que}{\begin{question}}
\newcommand{\cor}{\begin{corollary}}
\newcommand{\fac}{\begin{fact}}

\newcommand{\prf}{\noindent{\bf Proof:} }
\newcommand{\ethm}{\end{theorem}}
\newcommand{\elem}{\end{lemma}}
\newcommand{\epro}{\end{proposition}}
\newcommand{\edfn}{\bbox\end{definition}}
\newcommand{\erem}{\bbox\end{remark}}
\newcommand{\exam}{\bbox\end{example}}
\newcommand{\ecnj}{\bbox\end{conjecture}}
\newcommand{\eprb}{\bbox\end{problem}}
\newcommand{\eque}{\bbox\end{question}}
\newcommand{\ecor}{\end{corollary}}
\newcommand{\efac}{\end{fact}}
\newcommand{\eprf}{\bbox}
\newcommand{\beqn}{\begin{equation}}
\newcommand{\eeqn}{\end{equation}}
\newcommand{\wbox}{\mbox{$\sqcap$\llap{$\sqcup$}}}
\newcommand{\bbox}{\vrule height7pt width4pt depth1pt}
\newcommand{\qed}{\bbox}

\newcommand{\rarrow}{\rightarrow}
\newcommand{\larrow}{\leftarrow}
\newcommand{\grad}{\bigtriangledown}

\overfullrule=0pt
\def\setof#1{\lbrace #1 \rbrace}
\section{Introduction}
This paper is concerned with multivariate polynomials and entire functions with nonnegative
real coefficients.(All Taylor's series in this paper are taken at zero.) We continue the research, initiated in the recent papers \cite{newhyp}, \cite{my-vdw}, \cite{stoc-06}, \cite{fried}, \cite{gur} by the present author, on
``combinatorics and combinatorial applications
hidden in certain homogeneous polynomials with non-negative coefficients.''
Essentially, the main goal here is understanding how far one can push the approach from the above mentioned papers.
The following definition introduces the main notation of the paper.

\dfn
\begin{enumerate}
\item We denote by $Sim(n)$ the standard simplex in $R^n$:
$$
Sim(n) = \{(a_1,...,a_n): a_i \geq 0, 1 \leq i \leq n; \sum_{1 \leq i \leq n} a_i = 1.
$$

\item
We denote by $Pol_{+}(m,n)$ the convex cone of polynomials with nonnegative coefficients in $m$ variables
of total degree $n$;  the corresponding convex cone of homogeneous polynomials is denoted
as $Hom_{+}(m,n)$.\\
We denote by $Ent_{+}(m)$ the convex cone of entire functions on $C^m$ with nonnegative
Taylor's series.
\item
An entire function $f \in Ent_{+}(m)$ is called {\bf Strongly Log-Concave}
 if for all integer vectors $(c_1,...,c_m) \in Z_{+}^{m}$ the function
$(\partial x_1)^{c_1}...(\partial x_m)^{c_m} f$ is either zero or
$\log((\partial x_1)^{c_1}...(\partial x_m)^{c_m} f)$ is concave on $R_{+}^{m}$.
A set of {\bf Strongly Log-Concave} polynomials $p \in Pol_{+}(m,n)$ is denoted as $SLC(m,n)$ and
a set of {\bf Strongly Log-Concave} entire functions $f \in Ent_{+}(m)$
is denoted as $SLC(m)$.\\
\item
A (discrete) subset $S \subset Z^{m}$ is called $D$-convex if
$$
Conv(S) \cap  Z^{m} = S,
$$
where $Conv(S)$ is the convex hull of $S$ and $Z^m$ is the $m$-dimensional
integer lattice.\\
A map $G : Z^m \rightarrow [-\infty, +\infty]$ is called
$D$-concave if
$$
G(\sum_{1 \leq i \leq k < \infty} a_i Y_i) \geq \sum_{1 \leq i \leq k < \infty} a_i G(Y_i) 
$$
for all sequences $(a_1,...,a_k) \in Sym(k)$ and all vectors $Y_1,...,Y_k \in Z^m$ such that $\sum_{1 \leq i \leq k < \infty} a_i Y_i \in Z^m$.

{\it Our notion of $D$-convexity coincides with the notion of {\bf pseudo-convexity} from \cite{danilov}. As the term ``pseudo-convex" is already
occupied in the complex analysis, we think that the term {\bf $D$-convexity} is more appropriate (and informative).}

\item
The support of an entire function
\beqn \label{tayl}
f(x_1,...,x_m) = \sum_{(r_1,...,r_m) \in Z_{+}^m} a_{r_1,...,r_m} \prod_{1 \leq i \leq m } x_i^{r_{i}}
\eeqn
is defined as $supp(f) = \{(r_1,...,r_m) : a_{r_1,...,r_m} \neq 0 \}$.
\item
For an entire function $f \in Ent_{+}(m)$ and an integer vector $R = (r_1,...,r_m) \in Z_{+}^m$
we define $Der_{f}(R) =: (\partial x_1)^{r_1}...(\partial x_m)^{r_m} f(0)$.\\
In the notation of (\ref{tayl}), $Der_{f}(R) = a_{r_1,...,r_m} \prod_{1 \leq i \leq m} r_{i}!$

\end{enumerate}
\edfn
\xam
\begin{enumerate}
\item

First, we note that a homogeneous polynomial $p \in Hom_{+}(m,n)$ is log-concave on $R_{+}^{m}$ if and only if
the function $p^{\frac{1}{n}}$ is concave on $R_{+}^{m}$.\\
\item
A natural class of {\bf Strongly Log-Concave} homogeneous polynomials in $Hom_{+}(m,n)$ consists of {\bf H-Stable} polynomials:
a polynomial $p \in Hom_{C}(m,n)$ is called {\bf H-Stable} if $p(Z) \neq 0$ provided
$Re(Z) > 0$. It is easy to show and is well known that if $p \in Hom_{C}(m,n)$ is {\bf H-Stable} then
the polynomial $\frac{p}{p(x_1,...,x_m)} \in Hom_{+}(m,n)$ for any positive real vector $(x_1,...,x_m)$ and
$(\partial x_i)p$ is either zero or {\bf H-Stable}.
Consider an univariate polynomial $R(t) = \sum_{0 \leq i \leq k} a_i t^i, a_k \neq 0$ and the associated
homogeneous polynomial $p \in Hom_{+}(2,n), p(x,y) = \sum_{0 \leq i \leq k} a_i x^i y^{n-i}$.\\
Then $p$ is {\bf H-Stable} iff the roots of $R$ are non-positive real numbers, which shows that {\bf H-Stable} polynomials
are {\bf Strongly Log-Concave}.
\item Another, different from {\bf H-Stable}, class of {\bf Strongly Log-Concave} homogeneous polynomials in $Hom_{+}(m,n)$ consists of
Minkowski polynomials $Vol_{n}(\sum_{1 \leq i \leq m} x_i K_i)$, where $Vol_n$ stands for the standard volume
in $R^n$ and $K_1,...,K_m$ are convex compact subsets of $R^n$. The {\bf Strong Log-Concavity} of Minkowski polynomials
is essentially equivalent to the famous {\bf Alexandrov-Fenchel} inequalities \cite{Al38} for the mixed volumes.

\rem {\bf H-Stable} and Minkowski polynomials satisfy a seemingly stronger property:
they are invariant respect to the changes of variables $Y = AX$, where $A$ is a rectangular matrix with {\bf non-negative}
entries and without zero rows. We don't know whether such invariance holds in the general  {\bf Strongly Log-Concave} case.
\erem

\end{enumerate}
\exam

We are interested in the following natural question: {\it when the support $supp(f)$ of an
entire function $f \in Ent_{+}(m)$ is $D$-convex?}. Clearly, $supp(f)$ is $D$-convex if, for instance,
the map $\log(Der_{f}):Z^m \rightarrow [-\infty, +\infty)$ is $D$-concave. This is the case
for $f(x,y) = \sum_{1 \leq i \leq n} a_i x^i y^{n-i}, f \in Hom_{+}(2,n)$ such that the univariate polynomial
$R(t) = \sum_{1 \leq i \leq n} a_i t^i$ has only real roots. Spelling out the definition of $D$-concavity
gives us a reformulation of the famous Newton's inequalities.\\ 
In the case of {\bf Strongly Log-Concave} multivariate entire functions, the map
$\log(Der_{f})$ is not necessary $D$-concave.\\
We introduce the following map
\beqn \label{cap-gen}
C_{f}(r_1,...,r_m) = \inf_{x_i > 0} \frac{f(x_1,...,x_m)}{ \prod_{1 \leq i \leq m }(\frac{x_i}{r_i})^{r_{i}} } , (r_1,...,r_m) \in Z_{+}^m
\eeqn
It is easy to show that if $f \in Ent_{+}(m)$ and $\log(f)$ is concave on $R_{++}^m$ then $\log(C_{f})$ is $D$-concave.\\
Therefore, the $D$-convexity of the support $supp(f)$ would follow from the property
\beqn \label{posit}
C_{f}(R) > 0 \Leftrightarrow Der_{f}(R) > 0.
\eeqn
We prove in this paper a sharp quantative version of (\ref{posit}):
\beqn \label{exp-waer}
\prod_{1 \leq i \leq m} \frac{r_i !}{r_{i}^{r_{i}}} C_{f}(R) \geq Der_{f}(R) \geq \exp \left(-(\sum_{1 \leq i \leq m} r_i) \right) C_{f}(R).
\eeqn
The inequalities (\ref{exp-waer}) (and their more refined versions) generalize the {\bf Van Der Waerden-Falikman-Egorychev} lower bound on
the permanent of doubly-stochastic matrices \cite{fal}, \cite{ego} and used in this paper to prove Newton-like inequalities for
{\bf Strongly Log-Concave} entire functions.

\section{Univariate Newton-like Inequalities}

\subsection{Propagatable sequences (weights)}
\dfn
Let us define the following closed subset of $R^{n+1}$ of log-concave sequences:
$$
LC = \{(d_0,...,d_n) : d_i \geq 0, 0 \leq i \leq n; d_i^{2} \geq d_{i-1}d_{i+1}, 1 \leq i \leq n-1 \}.
$$
We also associate with a given positive vector $(c_0,...,c_n)$ the weighted shift operator $Shift_{\bf{c}} : R^{n+1} \longrightarrow R^{n+1}$,
$$
Shift_{\bf{c}}((x_0,...,x_n)^{T}) = (c_{0}x_1,...,c_{n-1}x_{n},0)^{T}.
$$
If $c$ is the vector of all ones,
then $Shift_{\bf{c}} =: Shift$.\\
A positive finite sequence $(b_0,...,b_n)$ is called {\bf propagatable} if the following implication
holds:\\
\beqn \label{fin-prop}
(p^{(0)}(0)b_0,...,p^{(n)}(0)b_n ) \in LC \Longrightarrow (p^{(0)}(t)b_0,...,p^{(n)}(t)b_n ) \in LC , t \geq 0,
\eeqn
where $p$ is a polynomial of degree at most $n$.
\edfn
Analogously, we define infinite {\bf propagatable} sequences by considering infinite log-concave sequences and entire functions in (\ref{fin-prop}).

\pro \label{ppg}
Let $c_0,...,c_{n-1}$ be a nonnegative sequence. Then $exp(t(Shift_{\bf{c}}))(LC) \subset LC$ for all $t \geq 0$ if and only if
$$
2c_{i} \geq c_{i+1} + c_{i-1}, 1 \leq i \leq n-2; 2c_{n-1} \geq c_{n-2}.
$$
(In other words, the infinite sequence $(c_0,...,c_{n-1},0,...)$ is concave).
\epro

\prf
\begin{enumerate}
\item {\it The "only if" part:} Consider the linear system of differential equations :
$$
X^{\prime}(t) = Shift_{\bf{c}} X(t): X(0) = (1,1,...,1), X(t) = (X_{0}(t),...,X_{n}(t)).
$$
Suppose that $exp(tShift_{\bf{c}})(LC) \subset LC, t \geq 0$ , i.e $X(t) \in LC: t \geq 0$.\\
Define the following smooth functions: 
$$
r_{i}(t) = (X_{i}(t))^2 - X_{i+1}(t) X_{i-1}(t), 1 \leq i \leq n-1.
$$
It follows that $r_{i}(0) = 0$ and $r_{i}(t) \geq 0, t \geq 0$. Therefore
$r_{i}^{\prime}(0) \geq 0$. Thus
$$
0 \leq r_{i}^{\prime}(0) = 2c_{i} - c_{i+1} - c_{i-1}, 1 \leq i \leq n-2; 0 \leq  r_{n-1}^{\prime}(0) = 2c_{n-1} - c_{n-2}.
$$
\item {\it The "if" part:} As $exp(A) = \lim_{n \rightarrow \infty} (I + \frac{A}{n})^n$,
thus it is sufficient to prove that\\
$(I + tShift_{\bf{c}})(LC) \subset LC$ for all $t \geq 0$,
which is done by straigthforward derivations.\\

\end{enumerate}
\eprf
\rem
The observation that $(I + Shift)(LC) \subset LC$ is probably well known; we have learned it from Julius Borcea.
\erem
\thm \label{prop}
Let $(b_0,...,b_k)$ be a positive sequence. Define $c_i = \frac{b_{i}}{b_{i+1}}, 0 \leq i \leq k-1$.
The sequence $(b_0,...,b_k)$ is {\bf propagatable} iff the infinite sequence $(c_0,...,c_{k-1},0,...)$ is concave.
\end{theorem}

\prf
Define a vector function $Mom_{b}(t) =(b_{0}p^{(0)}(t),...,b_{n}p^{(n)}(t))$.
Clearly, $Mom_{b}(t)$ solves the following system of  linear differential equations:
$$
Mom_{b}(t)^{\prime} = Shift_{{\bf c}}(Mom_{b}(t)).
$$
Therefore $(b_0,...,b_n)$ is {\bf propagatable} iff\\
$exp(t(Shift_{\bf{c}}))(LC) \subset LC$ for all $t \geq 0$. The result now follows from Proposition (\ref{ppg}).
\eprf

The following result follows fairly directly from Theorem (\ref{prop}).
\cor \label{infin}
Let $(b_0,...,b_k,...)$ be a positive infinite sequence. Define $c_i = \frac{b_{i}}{b_{i+1}}, 0 \leq i < \infty$.
The sequence $(b_0,...,b_k,...)$ is {\bf propagatable} iff the infinite sequence $(c_0,...,c_{k-1},...)$ is concave.
\ecor

\xam
A polynomial $p(t) = \sum_{0 \leq i \leq k} a_i t^i$ with nonnegative coefficients is called $n$-Newton for $n \geq k$
if
\beqn \label{prav}
d_i^{2} \geq d_{i-1} d_{i+1} : 1 \leq i \leq k-1, d_i =: \frac{a_i}{{n \choose i}}.
\eeqn
Or, in other words, the vector $(p^{(0)}(0)b_0,...,p^{(k)}(0)b_k ) \in LC$,
where $b_{i} = (n-i)!$.\\
As $c_i = \frac{b_{i}}{b_{i+1}} = n-i : 0 \leq i \leq k-1$ hence it follows from Theorem (\ref{prop})
that\\
 $(p^{(0)}(t)b_0,...,p^{(k)}(t)b_k ) \in LC : t \geq 0$. Equivalently,
\beqn \label{prim}
(p^{(i+1)}(t))^2 \geq \frac{n-i}{n-i -1} p{(i)}(t)p^{(i+2)}(t) : t \geq 0, i \leq k-2,
\eeqn
which means
that the functions $\sqrt[n-i]{p^{(i)}}: 0 \leq i \leq k$ are concave on $R_{+}$.\\

Let $f \in Ent_{+}(1)$ be entire univariate function, $f(t) = \sum_{0 \leq i < \infty} a_i t^i$. \\
A natural generalization of the $n$-Newton property,
i.e. when $n \rightarrow \infty$, is the log-concavity of
the infinite sequence $f^{(0)}(0),...,f^{(k)}(0),...$. Corollary (\ref{infin}) proves
that this property is equivalent to {\bf Strong Log-Concavity} of $f$.\\
We collect the above observations in the following proposition.

\pro \label{alter}

\begin{enumerate}
\item
A polynomial $p$ with nonnegative coefficients is $n$-Newton, where $n \geq deg(p)$, iff
the functions $\sqrt[n-i]{p^{(i)}}: 0 \leq i \leq k$ are concave on $R_{+}$.\\

Let us $n$-homogenize the univariate polynomial $p$, i.e. put $R(x,y) = y^n p(\frac{x}{y})$.
Then, $R \in Hom_{+}(2,n)$ and the functions $\sqrt[n-i]{p^{(i)}}: 0 \leq i \leq k$ are concave on $R_{+}$
if and only if the polynomial $R$ is {\bf Strongly Log-Concave}.

\item
An entire function $f \in Ent_{+}(1)$ is {\bf Strongly Log-Concave} iff
the infinite sequence $f^{(0)}(0),...,f^{(k)}(0),...$ is log-concave.
\end{enumerate}
\epro
\rem
The standard Newton Inequalities correspond to the case $n = deg(p)$ and hold if, for instance, the roots of $p$ are real.
It was proved by G. C. Shephard in \cite{shep} that a polynomial $p$ is $n$-Newton iff
$p(t) = Vol_{n}(t K_1 + K_2)$ for some convex compact subsets(simplices) $K_1, K_2 \subset R^n$. This remarkable
result can be used (see \cite{mixvol} and \cite{lig}) for alternative short proofs of Proposition (\ref{alter}) and
Liggett's
convolution theorem, which states that $pq$ is $m+n$-Newton provided that $p$ is $n$-Newton
and $p$ is $m$-Newton.\\
The literature on univariate Newton Inequalities is vast, we refer the reader to the recent survey \cite{roma}.
But the results presented here seem to be new, nothing of the kind is mentioned in \cite{roma}.\erem
\exam 
\section{Multivariate Case}
The main upshot of Proposition(\ref{alter}) is that in the univariate case as well in the
bivariate homogeneous case the following equivalence holds:\\

``$f$ is {\bf Strongly Log-Concave}" $\Longleftrightarrow$ ``the map $\log(Der_f)$ is $D$-concave".\\

In the general multivariate case both implication fail.
\xam \label{cycle}

\begin{enumerate}
\item
Consider the polynomial
$p(x_1,...,x_{2n}) = (x_1 + x_2)(x_2 + x_3)...(x_{2n-1} + x_{2n})(x_{2n} + x_1)$. Clearly, it is {\bf H-Stable}.
 Consider three vectors:
$R_0 = (1,...,1), R_1 =(2,0,2,...,0,2), R_2 = (0,2,...,0,2); 2 R_0 = R_1 + R_2$. By direct inspection,
$Der_{p}(R_0) = 2, Der_{p}(R_1) = Der_{p}(R_2) = 2^n$. Which gives
\beqn \label{decr}
\log(Der_p(\frac{1}{2}(R_1 + R_2))) = \frac{1}{2} \left(\log(Der_p(R_1)) + \log(Der_p(R_2)) \right) - (n-1)  \log(2).
\eeqn
\item {\it Alexandrov-Fenchel Inequalities}.
\\

Consider a homogeneous {\bf Strongly Log-Concave} polynomial $p \in Hom_{+}(m,n)$ and fix
a non-negative integer vector $R = (r_1,r_2,...,r_m), \sum_{1 \leq i \leq m} r_i = m$.
Define the following  polynomial $q \in Hom_{+}(2,n- \sum_{3 \leq i \leq m} r_i)$,
$$
q(x_1,x_2) = (\partial x_3)^{r_3}...(\partial x_m)^{r_m} p(x_1,x_2,0,...,0).
$$
Then $q$ is either zero or {\bf Strongly Log-Concave}. This observation leads to the following
inequalities: if both vectors 
$$
R_1 =(r_1 +1, r_2 -1, r_3,...,r_m), R_2 = (r_1 -1, r_2 +1, r_3,...,r_m)
$$ 
are non-negative then
\beqn \label{min-con}
Der_{p}(R) = Der_{p}(\frac{1}{2}(R_1 + R_2)) \geq (Der_{p}(R_{1}))^{\frac{1}{2}} (Der_{p}(R_{2}))^{\frac{1}{2}}
\eeqn

\item
Consider $p \in Hom_{+}(4,4), p(x_1,x_2,x_3,x_4) = x_1 x_2 x_3 x_4 + \frac{1}{4} ((x_1 x_2)^2 + (x_3 x_4)^2)$.
Here the map $\log(Der_f)$ is $D$-concave but the polynomial $p$ is not log-concave on $R_{+}^4$.

\end{enumerate}
\exam

We prove in this paper that in the general multivariate case if $f$ is {\bf Strongly Log-Concave} then
the map $\log(Der_f)$ is ``almost" $D$-concave.

\subsection{Generalized Van Der Waerden-Falikman-Egorychev lower bounds}
This section follows the recent inductive approach by the author \cite{my-vdw}.
\dfn
For an entire function $f \in Ent_{+}(n)$ we define its {\bf Capacity} as
\beqn \label{cap}
Cap(f) =  \inf_{x_i > 0} \frac{p(x_1,...,x_n)}{ \prod_{1 \leq i \leq n }x_i }
\eeqn
\edfn
We need the following elementary result:
\lem \label{dur}
Consider a function $f:R_{+} \rightarrow R_{+}$
such that the derivative $f^{\prime}(0)$ exists.
\begin{enumerate}
\item If $f^{\frac{1}{k}}$ is concave on $R_{+}$  for $k > 1$ then $f^{\prime}(0) \geq (\frac{k-1}{k})^{k-1} \inf_{t>0} \frac{f(t)}{t}$.\\
\item If $f$ is log-concave on $R_{+}$  then $f^{\prime}(0) \geq \frac{1}{e} \inf_{t>0} \frac{f(t)}{t}$.\\
If, additionally, the function $f$ is analytic and $f^{\prime}(0) = \frac{1}{e} \inf_{t>0} \frac{f(t)}{t}$ then
$f(t) = exp(at), a > 0$.

\item Let $R(t) =a_0+...+a_n t^n$ be a strongly log-concave on $R_{+}$ univariate polynomial with nonnegative coefficients:\\
$ G(i)^2 \geq G(i-1) G(i+1): 1 \leq i \leq n-1, G(i) = a_i i!$.\\
Then $f^{\prime}(0) \geq L(n) \inf_{t>0} \frac{f(t)}{t}$, where $L(n) = (\inf_{t>0} \frac{exp_{n}(t)}{t})^{-1}$
and the truncated exponential is defined as $exp_{n}(t) = 1+...+ \frac{1}{n!} t^n$. (Note that $exp_{n}$ is strongly log-concave on $R_{+}$.)
\end{enumerate}
\elem
\prf
\begin{enumerate}
\item If $f(0) = 0$ then, obviously, $f^{\prime}(0) \geq \inf_{t>0} \frac{f(t)}{t}$. Therefore, we can assume that $f(0) = 1$.
As $f^{\frac{1}{k}}$ is concave and non-negative on $R_{+}$ thus\\
$f(t) \leq (1 + \frac{f^{\prime}(0)}{k} t)^k, t \geq 0$.\\
The standard calculus gives us for $l(t) = (1 + \frac{f^{\prime}(0)}{k} t)^k$ that\\
$$
\inf_{t > 0}\frac{l(t)}{t} = f^{\prime}(0) (g(k))^{-1}, g(k) = \left( \frac{k-1}{k} \right)^{k-1}.
$$
As $\inf_{t>0} \frac{f(t)}{t} \leq \inf_{t>0} \frac{l(t)}{t}$, we deduce that $f^{\prime}(0) \geq g(k) \inf_{t>0} \frac{f(t)}{t}$.
\item
As in the proof above, we can assume that $f(0) =1$. It follows from the log-concavity
that $f(t) \leq exp(f^{\prime}(0) t), t \geq 0$. It is easy to see that
$$
\inf_{t>0} \frac{f(t)}{t} \leq \inf_{t>0} \frac{exp(f^{\prime}(0) t)}{t} = f^{\prime}(0) exp(1) = \frac{exp(f^{\prime}(0) s)}{s}, s = (f^{\prime}(0))^{-1}.
$$
Therefore, $f^{\prime}(0) \geq \frac{1}{e} \inf_{t>0} \frac{f(t)}{t}$.\\
If $f^{\prime}(0) = \frac{1}{e} \inf_{t>0} \frac{f(t)}{t}$ then, using the log-concavity again, we get
that $f(t) = exp(f^{\prime}(0) t), 0 \leq t \leq s$. If $f$ is analytic then $f(z) = exp(az), z \in C, a = f^{\prime}(0) > 0$.
\item Again, assume WLOG that $R(0) = 1$. It follows then from the strong log-concavity that
$$
R(t) \leq  1+...+ \frac{1}{n!} t^n = exp_{n}(t), t \geq 0.
$$
The rest of the proof is now as above.
\end{enumerate}

\eprf

\cor \label{capas}
Let $f \in Ent_{+}(n+1)$ and $g_{n}(x_1,...,x_n) = (\partial x_{n+1})p(x_1,...,x_n,0)$.\\
If $f$ is log-concave on $R_{+}^{n+1}$ then
\beqn \label{exp-bnd}
Cap(q_{n}) \geq \frac{1}{e} Cap(f).
\eeqn
If $p \in Hom_{+}(n+1,n+1)$ is log-concave on $R_{+}^{n+1}$ then

\beqn \label{g-fun}
Cap(q_{n}) \geq g(n+1) Cap(p), \ \mbox{where} \ g(k)=:\left( \frac{k-1}{k} \right)^{k-1}.
\eeqn
\ecor
\prf
We need to prove that $(\partial x_{n+1})p(x_1,...,x_n,0) \geq \frac{1}{e} Cap(p) \prod_{1 \leq i \leq n} x_i $.
Define an univariate log-concave entire function
$R(t) = f(x_1,...,x_n,t)$.\\
Then
$R(t) \geq Cap(p) t \prod_{1 \leq i \leq n} x_i: t \geq 0$ and $R^{\prime}(0) = (\partial x_{n+1})f(x_1,...,x_n,0)$.\\
It follows from the second item in Lemma(\ref{dur}) that
$$
(\partial x_{n+1})p(x_1,...,x_n,0) \geq \frac{1}{e} Cap(p) \prod_{1 \leq i \leq n} x_i.
$$
The inequality (\ref{g-fun}) is proved in the very same way, using the first item in Lemma (\ref{dur}) and the fact that if $p \in Hom_{+}(n+1,n+1)$ is
log-concave on $R_{+}^{n+1}$
then also $p^{\frac{1}{n+1}}$ is concave on $R_{+}^{n+1}$.
\eprf
We use below the following notation:
$$
vdw(n) = \frac{n!}{n^n}.
$$
\thm \label{main-thm}
\begin{enumerate}
\item
Let $f \in Ent_{+}(n)$ be {\bf Strongly Log-Concave} entire function in $n$ variables.
Then the following inequality holds:
\beqn \label{ent-vdw}
Cap(f) \geq \frac{\partial^n}{\partial x_1...\partial x_n} f(0) \geq \frac{1}{e^n} Cap(f)
\eeqn
{\it Note that the right inequality in (\ref{ent-vdw}) becomes equality if
$f = exp(\sum_{1 \leq i \leq n} a_i x_i)$ where $a_i > 0, 1 \leq i \leq n$.}
\item
Let a homogeneous polynomial $p \in Hom_{+}(n,n)$ be {\bf Strongly Log-Concave}.
Then the next inequality holds:
\beqn \label{pol-vdw}
Cap(f) \geq \frac{\partial^n}{\partial x_1...\partial x_n} f(0) \geq vdw(n) Cap(p)
\eeqn
{\it Note that the right inequality in (\ref{pol-vdw}) becomes equality if
$p = (\sum_{1 \leq i \leq n} a_i x_i)^{n}$ where $a_i > 0, 1 \leq i \leq n$.}
\item
Let a polynomial $p \in Pol_{+}(n,n)$ be {\bf Strongly Log-Concave}.
Then the next inequality holds:
\beqn \label{pol-vdw1}
Cap(f) \geq \frac{\partial^n}{\partial x_1...\partial x_n} f(0) \geq \prod_{1 \leq  i \leq n} L(i) Cap(p),
\eeqn
where $L(n) = (\inf_{t>0} \frac{exp_{n}(t)}{t})^{-1}$.\\
({\it Note that $L(1)=1, L(2) = (1 + \sqrt{2})^{-1}$ and $L(n) > e^{-1}, n \geq 1$.})
\end{enumerate}
\ethm
\prf
\begin{enumerate}
\item
Define the following entire functions $q_{i} \in Ent_{+}(i)$:\\
$q_{n} = f$, $q_{i}(x_1,...,x_i) = \frac{\partial^{n-i}}{\partial x_{i+1}...\partial x_n}f(x_1,...,x_i,0,...,0)$.
Notice that $q_{1}^{\prime}(0) = \frac{\partial^n}{\partial x_1...\partial x_n} f(0) $.\\
By the definition of {\bf Strongly Log-Concavity}, these entire functions are either log-concave or zero.
Using the inequality (\ref{exp-bnd}), we get that
$$
Cap(q_{i}) \geq \frac{1}{e} Cap(q_{i+1}), 1 \leq i \leq n-1.
$$
Therefore
$$
\inf_{t>0} \frac{q_{1}(t)}{t} = Cap(q_{1}) \geq  (\frac{1}{e})^{n-1} Cap(f).
$$
Finally, using Lemma (\ref{dur}), we get that
$$
\frac{\partial^n}{\partial x_1...\partial x_n} f(0) = q_{1}^{\prime}(0) \geq \frac{1}{e} \inf_{t>0} \frac{q_{1}(t)}{t} \geq \frac{1}{e^n} Cap(f).
$$
\item
If a homogeneous polynomial $p \in Hom_{+}(n,n)$ is {\bf Strongly Log-Concave} then the polynomials
$q_{i} \in Hom_{+}(i,i)$, $\frac{\partial^n}{\partial x_1...\partial x_n} p(0) = Cap(q_{1})$ and $(q_{i})^{\frac{1}{i}}$ is concave on $R^i_{+}, 1 \leq i \leq n$.
It follows from the inequality (\ref{g-fun}) that
$$
\frac{\partial^n}{\partial x_1...\partial x_n} p(0) = Cap(q_{1}) \geq \prod_{2 \leq k \leq n} g(k) Cap(p) = \frac{n!}{n^n} Cap(p).
$$
\end{enumerate}

\eprf

\subsection{General monomials}

Consider an entire function $f \in Ent_{+}(m)$ and an integer non-negative vector
$R = (r_1,...,r_m)$ .\\
Assume WLOG that $R = (r_1,...,r_k,0,...,0) : r_i > 0, 1 \leq i \leq k; k \leq n$. Let us define
the entire function $f_{R} \in Ent_{+}(|R|_{1})$, where $|R|_{1} = r_1+...+r_k$.
$$
f_{(R)}(y_1,...,y_{|R|_{1}}) = f(e_{1}(y_1+...+y_{r_{1}})+...+e_{k}(y_{r_1+...+r_{k-1} +1} +...+ y_{r_1+...+r_{k}})),
$$
where $\{e_1,...,e_m\}$ is the standard basis in $C^m$. The following identity is obvious:
$$
(\partial x_1)^{r_1}...(\partial x_m)^{r_m} f(0) = (\partial y_1) ...(\partial y_{|R|_{1}}) f_{(R)} (0).
$$
{\bf Note that if the original entire function (homogeneous polynomial) $f$ is {\bf Strongly Log-Concave} ({\bf H-Stable}) then the same holds
for the entire function (homogeneous polynomial) $f_{(R)}$.}

It easily follows from the arithmetic-geometric means inequality that
\beqn \label{caps-con}
Cap(f_{(R)}) = C_{f}(r_1,...,r_m) =: \inf_{x_i > 0} \frac{f(x_1,...,x_m)}{ \prod_{1 \leq i \leq m }(\frac{x_i}{r_i})^{r_{i}} }
\eeqn

As we deal only with entire functions with the non-negative coefficients hence the following
inequality holds:
\beqn \label{upper-bn}
\left(\prod_{1 \leq i \leq m} vdw(r_{i})\right)C_{f}(r_1,...,r_m)  \geq (\partial x_1)^{r_1}...(\partial x_m)^{r_m} f(0)
\eeqn

Putting these observations together, we get the Corollary to Theorem(\ref{main-thm}).
\cor
\begin{enumerate}
\item
Let $f \in Ent_{+}(m)$ be {\bf Strongly Log-Concave} entire function in $m$ variables. Then for all
integer vectors  $R = (r_1,...,r_m) \in Z_{+}^{m}$ the next inequalities hold:
\beqn \label{mon-vdw}
\left(\prod_{1 \leq i \leq m} vdw(r_{i})\right)C_{f}(r_1,...,r_m)  \geq (\partial x_1)^{r_1}...(\partial x_m)^{r_m} f(0) \geq exp(-|R|_{1}) C_{f}(r_1,...,r_m)
\eeqn
\item
Let a homogeneous polynomial $p \in Hom_{+}(m,n)$ be {\bf Strongly Log-Concave}. Then for all
integer vectors $R = (r_1,...,r_m) \in Z_{+}^{m}, \sum_{1 \leq i \leq m} r_i = n$ the next inequalities hold:
\beqn \label{segod}
\left(\prod_{1 \leq i \leq m} vdw(r_{i})\right)C_{p}(r_1,...,r_m)  \geq (\partial x_1)^{r_1}...(\partial x_m)^{r_m} p(0) \geq vdw(n) C_{p}(r_1,...,r_m)
\eeqn
\end{enumerate}

\ecor

Let us recall the generalized Schrijver's inequality from \cite{my-vdw}.
\thm \label{schr}
Let $p \in Hom_{+}(n,n)$ be {\bf H-Stable}. Let us denote the degree of variable $x_i$ in the polynomial $p$ as $deg_{p}(i)$.\\
If $deg_{p}(i) \leq k \leq n, 1 \leq i \leq n$,
Then the next inequality holds:
\beqn
Cap(p) \geq \frac{\partial^n}{\partial x_1...\partial x_n} p(0) \geq (\frac{k-1}{k})^{(k-1)(n-k)} vdw(k) Cap(p)
\eeqn
\ethm
Combining Theorem (\ref{schr}) and observations (\ref{caps-con}), (\ref{upper-bn}) we get the following Corollary.

\cor

Let $p \in Hom_{+}(n,n)$ be {\bf H-Stable}. Assume that the degree of variable $x_i$ in the polynomial $p$, $deg_{p}(i) \leq k \leq n, 1 \leq i \leq n$.
Then the following inequalities hold:
\beqn
\left(\prod_{1 \leq i \leq m} vdw(r_{i})\right)C_{p}(r_1,...,r_m)  \geq (\partial x_1)^{r_1}...(\partial x_m)^{r_m} p(0) \geq  (\frac{k-1}{k})^{(k-1)(n-k)} vdw(k)C_{p}(r_1,...,r_m)
\eeqn
\ecor

\subsection{A lower bound on the inner products of {\bf H-Stable} polynomials}
\thm \label{int-supp}
Let us consider two {\bf H-Stable} polynomials $p, q \in Hom_{+}(m,n)$:
$$
p(x_1,...,x_m) = \sum_{r_1+...+r_m = n } a_{r_1,...,r_m} \prod_{1 \leq i \leq m } x_i^{r_{i}}, q(x_1,...,x_m) = \sum_{r_1+...+r_m = n} b_{r_1,...,r_m} \prod_{1 \leq i \leq m } x_i^{r_{i}},
$$
and a nonnegative vector $(l_1,...,l_m)$ such that $\sum_{1 \leq i \leq m} l_i = n$.

Let us assume that
\beqn \label{interse}
\inf_{x_i > 0, 1 \leq i \leq m} \frac{p(x_1,...,x_m)}{\prod_{ 1 \leq i \leq m} x_i^{l_{i}}} =: A > 0, 
\inf_{x_i > 0, 1 \leq i \leq m} \frac{q(x_1,...,x_m)}{\prod_{ 1 \leq i \leq m} x_i^{l_{i}}} =: B > 0.
\eeqn
Then the following inequality holds:
\beqn \label{new-war}
<p,g> =: \sum_{r_1+...+r_m = n } a_{r_1,...,r_m} b_{r_1,...,r_m} \geq AB \frac{vdw(nm)}{vdw(n)^{m}} 
\eeqn
\ethm

\prf
Let us consider a rational function $F = \prod_{1 \leq i \leq m} x_i^{n} p(x_1,...,x_m) q(\frac{1}{x_1},...,\frac{1}{x_m})$.
It is clear that, in fact, $F \in Hom_{+}(m,nm)$ and $F$ is {\bf H-Stable}. Note that
$$
(n!)^{m} \sum_{r_1+...+r_m = n } a_{r_1,...,r_m} b_{r_1,...,r_m} = (\partial x_1)^{n}...(\partial x_m)^{n} F(0).
$$
It follows from (\ref{interse}) that $C_{F}(n,...,n) \geq AB n^{nm}$. Using the right inequality in (\ref{segod}), we get
that
$$
 \sum_{r_1+...+r_m = n } a_{r_1,...,r_m} b_{r_1,...,r_m} \geq  AB \frac{vdw(nm)}{vdw(n)^{m}}.
$$
\eprf

\rem
\begin{enumerate}
\item
It is easy to see that the inequalities (\ref{interse}) holds for some vector $(l_1,....,l_m)$ if and only if the Newton polytopes, $Newt(p)$ and $Newt(q)$,
have non-empty intersection. (Recall that the Newton polytope $Newt(p)$ is the convex hull of the support $supp(p)$.)\\
One of the corollaries of Theorem (\ref{int-supp}) is the fact that the intersection $Newt(p) \cap Newt(q)$ is not empty iff
the intersection $supp(p) \cap supp(q)$ is not empty. There is alternative( and harder) way to prove this fact. It was proved in \cite{newhyp} and \cite{stoc-06} that
if $p$ is a {\bf H-Stable} polynomial then the Newton polytope $Newt(p)$ is the {\bf polymatroid}, based on some integer valued
submodular function.It follows from the celebrated Edmonds' result \cite{edm} that all the vertices of $Newt(p) \cap Newt(q)$ are integer.
Therefore, if $Newt(p) \cap Newt(q)$ is not empty then the exists an integer vector $(r_1,...,r_m) \in Newt(p) \cap Newt(q)$.
But all integer vectors in $Newt(p)$($Newt(q)$) belong to the support $supp(p)$($supp(q)$).\\
The inequality (\ref{new-war}) is unlikely sharp. We conjecture here a sharp version:
$$
\sum_{r_1+...+r_m = n } a_{r_1,...,r_m} b_{r_1,...,r_m} \prod_{1 \leq i \leq m} (r_i)! \geq  AB \frac{n!}{m^n}.
$$
\item
If {\bf H-Stable} polynomials $p,q \in Hom_{+}(m,n)$ are both multilinear, i.e. $deg_{p}(i), deg_{q}(i) \leq 1, 1 \leq i \leq m$,
then 
$$
<p,q> = (\partial x_1)...(\partial x_m)G(0),\ \mbox{where} \ G(x_1,...,x_m) = (\prod_{1 \leq i \leq m} x_i) p(x_1,...,x_m) q(\frac{1}{x_1},...,\frac{1}{x_m}).
$$
Note that the polynomial $G \in Hom_{+}(m,m)$ is {\bf H-Stable}, $Cap(G) \geq AB$ and $deg_{G}(i) \leq 2, 1 \leq i \leq m$. Using
Theorem (\ref{schr}), we get the following inequality:
\beqn \label{inner-pro}
<p,q> \ \geq  \ AB 2^{-m +1}
\eeqn 
The inequality (\ref{inner-pro}) is sharp for $m = 2n$.
\end{enumerate}
\erem

\section{Multivariate Newton Inequalities}
We start with the following simple fact.
\fac
If an entire function $f \in Ent_{+}(m)$ is log-concave on $R_{+}^m$ then the map $C_{f}$, defined as
$$
C_{f}(y_1,...,y_m) = \inf_{x_i > 0} \frac{f(x_1,...,x_m)}{\prod_{1 \leq i \leq m} (\frac{x_i}{y_i})^{y_{i}} }, y_i \geq 0
$$
is log-concave on $R_{+}^{m}$.
\efac
\prf
Assume WLOG that $y_i > 0, 1 \leq i \leq k \leq m$ and $y_j = 0, k+1 \leq j \leq m$. It follows from the monotonicity of $f$ that
$$
C_{f}(x_1,...,x_m) = \inf_{x_i > 0, 1 \leq i \leq k} \frac{f(x_1,...,x_k,0,...,0)}{\prod_{1 \leq i \leq k} (\frac{x_i}{y_i})^{y_{i}} }.
$$
Therefore $C_{f}(y_1,...,y_m) \geq a$ iff $\log(f(x_1 y_1,...,x_m y_m)) \geq \log(a) + \sum_{1 \leq i \leq m} y_i \log(x_i)$ for
all positive vectors $(x_1,...,x_m)$. The desired log-concavity follows now from the log-concavity of the function $f$ and of the logarithm.
\eprf

Let $Y =(r_1,...,r_m) \in Z_{+}^m$ be an integer vector. We use below the following notations:
$$
VDW(Y) = \prod_{1 \leq i \leq m } vdw(r_i), \ \mbox{where} \ vdw((r) = \frac{r!}{r^r}.
$$
\thm \label{main-ineq}
Let us consider integer vectors $Y_0, Y_1,...,Y_k \in Z_{+}^m$ such that
$$
Y_0 = \sum_{1 \leq i \leq k} a_i Y_i; a_i \geq 0,\sum_{1 \leq i \leq k} a_i =1.
$$

\begin{enumerate}
\item Suppose that the entire function $f \in Ent_{+}(m)$ is {\bf Strogly Log-Concave}. Then
\beqn \label{ent-new}
Der_{f}(Y_0) \geq \left(exp(-|Y_0|_{1}) \prod_{1 \leq i \leq k} (VDW(Y_i))^{-a_i} \right)     \prod_{1 \leq i \leq k} (Der_{f}(Y_i))^{a_i}
\eeqn
\item
If $p \in Hom_{+}(m,n)$ is {\bf Strogly Log-Concave} then
\beqn \label{hom-in}
Der_{f}(Y_0) \geq \left(vdw(n) \prod_{1 \leq i \leq k} (VDW(Y_i))^{-a_i} \right)    \prod_{1 \leq i \leq k} (Der_{f}(Y_i))^{a_i}
\eeqn
\item
If $p \in Hom_{+}(m,n)$ is {\bf H-Stable} and $deg_{p}(i) \leq k \leq n$ for all $1 \leq i \leq m$ then
\beqn \label{spar}
Der_{f}(Y_0) \geq \left( \left(\frac{k-1}{k} \right)^{(k-1)(n-k)} vdw(k) \prod_{1 \leq i \leq k} (VDW(Y_i))^{-a_i} \right)    \prod_{1 \leq i \leq k} (Der_{f}(Y_i))^{a_i}
\eeqn
\end{enumerate}

\ethm

\prf
We wiil prove only the inequality (\ref{ent-new}) as the other ones are proved in the same way.\\
Using the the right inequality in (\ref{mon-vdw}), we get that
$$
Der_{f}(Y_0) \geq exp(-|Y_0|_{1}) C_{f}(Y_0).
$$
Since the map $C_{f}$ is log-concave hence
$$
C_{f}(Y_0) \geq \prod_{1 \leq i \leq k} (C_{f}(Y_i))^{a_i}.
$$
Finally, we use the left inequality in (\ref{mon-vdw}):
$$
C_{f}(Y_i) \geq (VDW(Y_i))^{-1} Der_{f}(Y_i).
$$
\eprf

\cor \label{l-conv}
The support $supp(f)$ of {\bf Strogly Log-Concave} entire function $f \in Ent_{+}(m)$ is $D$-convex.
\ecor

\xam
\begin{enumerate}
\item
Let us consider the following vectors in $Z_{+}^n$:
$$
Y_0 = (1,1,...,1); Y_1 = (n,0,...,0),...,Y_n = (0,0,...,n).
$$
Note that $Y_0 = \sum_{1 \leq i \leq n} \frac{1}{n} Y_i$. If $p \in Hom_{+}(n,n)$ is {\bf Strogly Log-Concave} then
(\ref{hom-in}) gives the next inequality
$$
Der_{p}(Y_0) \geq   \prod_{1 \leq i \leq k} (Der_{f}(Y_i))^{\frac{1}{n}},
$$
which is attained on $p(x_1,...,x_n) = (x_1+...+x_n)^n$.\\
\item
Consider three vectors in $Z_{+}^{2n}$:\\
$$
Y_0 = (1,1,...,1); Y_1 = (2,...,2,0,...,0),Y_1 = (0,...,0,2,...,2), |Y_1|_1 = |Y_2|_{1} = 2n.
$$
If $p \in Hom_{+}(2n,2n)$ is {\bf H-Stable} and $deg_{p}(i) \leq 2 \leq 2n$ for all $1 \leq i \leq 2n$ then it follows from (\ref{spar}) that
\beqn \label{sharp}
Der_{p}(Y_0) \geq 2^{-n+1} \prod_{1 \leq i \leq 2} (Der_{p}(Y_i))^{\frac{1}{2}}.
\eeqn
The inequality (\ref{sharp}) is attained on the polynomial $p(x_1,...,x_{2n}) = (x_1 + x_2)(x_2 + x_3)...(x_{2n-1} + x_{2n})(x_{2n} + x_1)$.
\end{enumerate}

\exam

\section{Comments and Open problems}
\begin{enumerate}
\item The inequality (\ref{pol-vdw}) is a far going generalization of the famous Van der Waerden conjecture
on the permanent of doubly-stochastic matrices(\cite{minc}, \cite{fal}, \cite{ego} and the Bapat's conjecture \cite{bapat}),\cite{gur}.
See more on this combinatorial connection in \cite{stoc-06}, \cite{all}, \cite{fried}.\\
The Van der Waerden conjecture conjecture corresponds to {\bf H-Stable} polynomials
$$
Prod_{A}(x_1,...,x_n) = \prod_{1 \leq i \leq n} \prod_{1 \leq j \leq n} A(i,j) x_j,
$$
where $n \times n$ matrix is non-negative entry-wise and has no zero rows. If such a matrix is doubly-stochastic,
i.e. all its rows and columns sum to $1$, then $Cap(Prod_{A}) = 1$.

The convex relaxation approach to Newton-like inequalities in Theorem(\ref{main-ineq})
was introduced by the author in \cite{gur} for the determinantal polynomials $\det(\sum_{1 \leq i \leq m} x_i A_i)$,
where $A_1,...,A_m$ are $n \times n$ hermitian PSD matrices.The corresponding inequalities in \cite{gur} are weaker
than in the present paper.

\item Just the log-concavity of $f$ is not sufficient for $D$-convexity of its support $supp(f)$ even for univariate polynomials
with non-negative coefficients. Indeed, consider $p(t) = t + t^3$. The fourth root $\sqrt[4]{p(t)}$ is concave on $R_{+}$:\\
$$
(p^{(1)}(t))^3 - \frac{4}{3} p(t) p^{(2)}(t) = (1 + 3 t^2)^2 - \frac{4}{3} (t + t^3) 6t = (t^2 -1)^2 \geq 0.
$$
This example can be "lifted" to a ``bad" log-concave homogeneous polynomial $q \in Hom_{+}(4,4)$:\\
$q(x,y,v,w) = (x+y)^{3}(v+w) + (v+w)^{3}(x+y)$. It is easy to see that $Cap(q) = 2^5$ but 
$$
\frac{\partial^4}{\partial x \partial y \partial v \partial w} q(0) =0.
$$

\item In the case of {\bf H-Stable} polynomials, Corollary (\ref{l-conv}) can be made much more precise:\\
Define, for a subset $S \subset \{1,...,m\}$ and a polynomial $p \in Hom_{+}(m,n)$, the integer number $Deg_{p}(S)$
equal to the maximum total degree attained on variables in $S$.\\
Then the following relation holds:

\beqn \label{rado}
a_{r_1,...,r_m} > 0 \Longleftrightarrow \sum_{j \in S} r_j \leq Deg_{p}(S) : S \subset \{1,...,m\}, p \in Hom_{+}(m,n).
\eeqn
Additionaly, the integer valued map $Deg_{p}: 2^{\{1,...,m\}} \rightarrow \{0,...,n\}$ is {\bf submodular}.

The characterization (\ref{rado}), proved in \cite{newhyp}, is a far going generalization of the Hall-Rado theorems on the existence of perfect
matchings.\\
The paper \cite{stoc-06} provides algorithmic applications of this result:
strongly polynomial deterministic algorithms for the membership problem as for the support
as well for the Newton polytope of {\bf H-Stable} polynomials $p \in Hom_{+}(m,n)$, given as oracles.\\
We don't know whether (\ref{rado}) works for {\bf Strongly Log-Concave} homogeneous polynomials. But
it would follow from the following conjecture/question:
\cnj
Let $p \in Hom_{+}(3,n)$ be {\bf Strongly Log-Concave}. Then there exist convex compact subsets $K_1,K_2,K_3 \subset R^n$ such that
\beqn \label{min}
p(x_1,x_2,x_3) = Vol_{n}(x_1 K_1 + x_2 K_2 + x_3 K_3) : x_1,x_2,x_3 \geq 0
\eeqn
\ecnj
Or put more modestly:

\que
Which {\bf Strongly Log-Concave} polynomials $p \in Hom_{+}(3,n)$ allow the representation (\ref{min})?
\eque

{\it The Minkowski polynomials $Vol_n \in Hom_{+}(3,n), Vol_{n}(x_1 K_1 + x_2 K_2 + x_3 K_3)$ actually have seemingly stronger, than {\bf Strong Log-Concavity}, property:\\
the polynomials $\prod_{1 \leq j \leq r < n} (\sum_{1 \leq i \leq 3} (a_{i,j}\partial x_i) Vol_{n}$ are either zero or log-concave on $R_{+}^3$ provided that
$a_{i,j} \geq 0$.\\}

\item
\que Is the set of {\bf Strongly Log-Concave} entire functions closed under the multiplications?\\
If true it would imply that the Minkowski sum $supp(f) + supp(g) = supp(fg)$ is $D$-convex if
$f,g$ are {\bf Strongly Log-Concave}.
\eque
We note that the problem of $D$-convexity of Minkowski sums of $D$-convex sets
was studied from a combinatorial point of view in \cite{danilov}.

\item
\que
What are the asymptotically exact constants in Theorem(\ref{main-ineq})? Are the cyclic polynomials in Example(\ref{cycle})
extremal?
\eque
\item
Can recently refuted Okounkov's conjecture \cite{okoun}, in the representation theory, on log-concavity of multiplicities
be fixed/generalized in the way similar to Theorem(\ref{main-ineq})?

\item Stable multivariate polynomials form a backbone of
linear multivariate control. If $p \in HSP_{+}(n,n)$ then
$Cap(p) = \inf_{Re(z_{i}) > 0} \frac{|p(z_1,...,z_n)|}{\prod_{1 \leq i \leq n} Re(z_{i})}$.
In other words, the capacity can be viewed as a measure of stability. What is
a meaning of capacity if terms of control/dynamics or in terms of the corresponding hyperbolic PDE?
\item Can our results be reasonably generalized to the fractional derivatives?

\end{enumerate}

\end{document}